\documentclass[12pt]{article}
 \usepackage{a4, latexsym, amsfonts}
 \newtheorem{theorem}{Theorem}[section]
\newtheorem{lemma}{Lemma}[section]
 \newtheorem{definition}{Definition}[section]
 \newtheorem{corollary}{Corollary}[section]

\newcommand{\supp}{{\rm supp}}

 \newenvironment{proof}{\trivlist
      \item[\hskip\labelsep
      {\itshape Proof.}]\normalfont}
      {\hspace*{\fill}$\Box$\endtrivlist}

\begin{document}

\title{Subsequence Sums of  Zero-sum free Sequences II}
\author{ Pingzhi Yuan \\
{\small School of Mathematics, South China Normal University
  , Guangzhou 510631, P.R.CHINA}\\
 {\small e-mail  mcsypz@mail.sysu.edu.cn} }

\date{}
\maketitle
 \edef \tmp {\the \catcode`@}
   \catcode`@=11
   \def \@thefnmark {}

    \@footnotetext { Supported by NSF of China (No. 10571180).}

    \begin{abstract} Let $G$ be a finite  abelian group, and let $S$ be a sequence of elements in $G$. Let
    $f(S)$   denote the number of elements in $G$ which can be
    expressed as the sum over a nonempty subsequence of $S$. In this
    paper, we
    determine  all the sequences $S$ that  contains no zero-sum subsequences
     and  $f(S)\leq 2|S|-1$.


    \end{abstract}

\section{Introduction}

Let $G$ be a finite abelian group (written additively)throughout the
present paper. $\mathcal{F}(G)$ denotes the free abelian monoid with
basis $G$, the elements of which are called $sequences$ (over $G$).
A  sequence of not necessarily distinct elements from $G$ will be
written in the form $S=g_1\cdot\, \cdots \,\cdot g_k=\prod_{i=1}^k
g_i=\prod_{g\in G}g^{\mathsf v_g(S)}\in\mathcal{F}(G)$, where
$\mathsf v_g(S) \geq 0$ is called the $multiplicity$ of $g$ in $S$.
Denote by $|S|= k$ the number of elements in $S$ (or the $length$ of
$S$) and let $\supp(S) = \{g\in G: \, \mathsf v_g(S)>0\}$ be the
$support$ of $S$.

We say that $S$ contains some $g\in G$ if $\mathsf v_g(S)\geq1$ and
a sequence $T\in\mathcal{F}(G)$ is a $subsequence$ of $S$ if
$\mathsf v_g(T) \leq \mathsf v_g(S)$ for every $g\in  G$, denoted by
$T|S$. If $T|S$, then let $ST^{-1}$ denote the sequence obtained by
deleting the terms of $T$ from $S$. Furthermore, by $\sigma(S)$ we
denote the sum of $S$, (i.e. $\sigma(S) = \sum_{ i=1}^k g_i =
\sum_{g\in G}\mathsf v_g(S)g\in G$). By $\sum(S)$ we denote the set
consisting of all elements which can be expressed as a sum over a
nonempty subsequence of $S$, i.e.
$$\sum(S) =\{\sigma(T): T \,{\rm is\, a\, nonempty\, subsequence\, of\,} S\}.$$
We write $f(S)=|\sum(S)|$,  $<S>$  for the subgroup of $G$ generated
by all the elements of $S$.

 Let $S$ be a sequence over $G$. We call $S$ a $zero-sum$ $sequence$ if $\sigma(S) =
 0$, a $zero-sum \,free$ $sequence$ if  $\sigma(W)\neq 0$ for any subsequence $W$ of $S$,
  and $squarefree$ if $\mathsf v_g(S)\leq1$ for every $g\in G$. We
  denote by $\mathcal{A}^\star(G)$ the set of all zero-sum free
  sequences in $\mathcal{F}(G)$.

Let $D(G)$ be the Davenport's constant of $G$, i.e., the smallest
integer $d$ such that every sequence $S$ over $G$ with $|S|\geq d$
satisfies $0\in\sum(S)$.  For every positive integer $r$ in the
interval $\{1,\,  \ldots, \,D(G)-1\}$, let
\begin{equation} f_G(r)=\min_{S, \,|S|=r}f(S),
\end{equation} where $S$ runs over all zero-sum free sequences of $r$
elements in $G$. How does the function $f_G$ behave?

In 2006, Gao and Leader proved the following result.

{\bf Theorem A} {\bf \cite{GL06}} {\it Let $G$ be a finite abelian
group of exponent $m$. Then

(i) If $1\leq r\leq m-1$ then $f_G(r)=r$.

(ii) If $\gcd(6,\, m)=1$ and $G$ is not cyclic then $f_G(m)=2m-1$.}

Recently, Sun\cite{sun07} showed that $f_G(m)=2m-1$ still holds
without the  restriction that $\gcd(6, \, m)=1$.

  Using some techniques from the author \cite{yu07}, the author \cite{yu08}
proved the following two theorems.

{\bf Theorem B}\cite{yu08,OW77} {\it Let $S$ be a zero-sum free
sequence over $G$ such that $<S>$ is not a cyclic group, then
$f(S)\geq 2|S|-1$.}

{\bf Theorem C} \cite{yu08} {\it Let  $S$ be a zero-sum free
sequence over $G$ such that $<S>$ is not a cyclic group and  $f(S)=
2|S|-1$. Then $S$ is one of the following forms

(i) $S=a^x(a+g)^y, \, x\geq y\geq 1$, where $g$ is an element of
order 2.

(ii) $S=a^x(a+g)^yg,\,  x\geq y\geq1$,  where $g$ is an element of
order 2.

(iii) $S=a^xb, \,x\geq1$.}

However, Theorem B is an old theorem of Olson and White \cite{OW77}
which has been overlooked by the author. For more recent progress on
this topic, see \cite{GP08,Pi09,yu09}.

The main purpose of the present paper is to determine all the
sequences $S$  over a finite abelian group such that $S$ contains no
zero-sum subsequences     and  $f(S)\leq 2|S|-1$. To begin with, we
need the notation of $g$-smooth.

\begin{definition} {\rm [7, Definition 5.1.3]} A sequence $S\in \mathcal{F}(G)$ is called
$smooth$ if $S=(n_1g)(n_2g)\cdot \, \cdots\,\cdot (n_lg)$, where
$|S|\in \mathbb{N}, \, g\in G, \, 1=n_1\leq \cdots\leq n_l,
\,n=n_1+\cdots+n_l<\mbox{ord}(g)$ and $\sum(S)=\{g,  \ldots, \,ng\}$
( in this case we say more precisely that $S$ is $g$-smooth).
\end{definition}

We have
\begin{theorem} Let $G$ be a finite abelian group and let $S$ be a zero-sum free sequence over
$G$ with $f(S)\le  2|S|-1$. Then $S$ has one of the following forms:

(i) $S$ is $a$-smooth for some $a\in G$.

(ii) $S=a^kb$, where $k\in\mathbb{N}$ and $a,  b\in G$ are distinct.

(iii) $S=a^kb^l$,  where $k\ge l \ge 2$ and $a, b \in G$ are
distinct with $2a = 2b$.

(iv) $S=a^kb^l(a-b) $,  where $k\ge l \ge 2$ and $a, b \in G$ are
distinct with $2a = 2b$.
\end{theorem}

For a sequence $S$ over $G$ we call
$$\mathsf h(S) =\max\{\mathsf v_g(S)|  g \in  G\}\in [0, |S|]$$
$$ the\, maximum\, of\, the\, multiplicities\, of\, S.$$
Let $S = a^xb^yT$ with $x\ge y \ge \mathsf h(T)$, then Theorem
1.1(i) can be stated more precisely as that $S$ is $a$-smooth or
$b$-smooth.
\section{Some Lemmas}

Let $\emptyset\neq G_{0}
 \subseteq G$ be a subset of $G$ and $k\in \mathbb{N}$. Define
 $$\mathsf f(G_{0},\,k)=\min\{f(S): \,S\in \mathcal{F}(G_{0})\, \,{\rm ~
 zero-sum free,\,
 squarefree ~and ~}\,|S|=k \} $$  and set
 $\mathsf f(G_{0},\,k)=\infty$, if there are no sequences over $G_{0}$ of the
 above form.
\begin{lemma}

Let \ $G$ \ be a finite abelian group.

\begin{enumerate}

\smallskip

\item If $k \in \mathbb N$ and $S = S_1 \cdot\, \cdots \,\cdot S_k \in
      \mathcal{A}^\star(G)$, then

      \[
      f(S) \ge f(S_1) + \cdots + f(S_k) \,.
      \]

\smallskip

\item If \ $G_0 \subset G$, \ $k \in \mathbb N$ and \ $\mathsf
f(G_0,\,k)>0$, then

$$
\mathsf f(G_0,\,k) \ \left\{\begin{array}{ll}
=1\,, \quad &\mbox{if} \quad k=1\,,\\
=3\,, \quad &\mbox{if} \quad k=2\,,\\
\ge 5\,, \quad &\mbox{if} \quad k=3\,,\\
\ge 6\,, \quad &\mbox{if} \quad k=3 \quad \mbox{and} \quad 2g \ne 0 \quad \mbox{for all} \quad g \in G_0\,,\\
\ge 2k\,, \quad &\mbox{if} \quad k \ge 4\,.
\end{array}\right.
$$

\end{enumerate}

\end{lemma}

\begin{proof}

1.  See \cite[Theorem 5.3.1]{Ge-HK06a}.

2. See \cite[Corollary 5.3.4]{Ge-HK06a}.

\end{proof}

\begin{lemma}
Let $a, \, b$ be  two distinct elements in an abelian group $G$ such
that  $a^2b^2\in\mathcal{A}^\star(G), \, 2a\neq 2b, a\neq 2b$, and
$b\neq 2a$.  Then $f(a^2b^2)=8$.
      \end{lemma}

\begin{proof} It is easy to see that $a, \, 2a, \, b, \, 2b, \, a+b,
\, a+2b, \, 2a+b, \, 2a+2b$ are all the distinct elements in
$\sum(a^2b^2)$. We are done.
\end{proof}

\begin{lemma}  Let $S=a^kb$ be a zero-sum free sequence over $G$. If
$S=a^kb$ is not $a$-smooth, then $f(S)=2k+1$.

\end{lemma}
\begin{proof} The assertion follows from the fact that $a, \,\ldots, \,ka, \,b,  \,a+b, \,\ldots, \,ka+b$
are all the distinct elements in $\sum(a^kb)$.\end{proof}

\begin{lemma} {\rm [10, Lemma 4]} Let $S$ be a zero-sum free sequence over $G$. If there is
some element $g$ in $S$ with order $2$, then $f(S)\geq2|S|-1$.
\end{lemma}

\begin{lemma}  Let $k\geq l\geq2$ be two
integers, and let $a$ and $b$ be two distinct elements of $G$ such
that $a^kb^l\in \mathcal{A}^\star(G)$ and $a^kb^l$ is not smooth.
Then we have

(i) If $2a\neq 2b$, then $f(a^kb^l)\geq 2(k+l)$.

(ii) If $2a=2b$, then $f(a^kb^l)= 2(k+l)-1$.
\end{lemma}
\begin{proof} If $nb\neq sa$ for any $n$ and $s$ with $1\leq n\leq
l$ and $1\leq s\leq k$, then $ra+sb, \,r+s\neq 0, \,0\leq r\leq k,\,
0\leq s\leq b$ are all the distinct elements in $\sum(a^kb^l)$, and
so
$$f(a^kb^l)=kl+k+l\geq 2(k+l).$$
Now we assume that $nb= sa$ for some $n$ and $s$ with $1\leq n\leq
l$ and $1\leq s\leq k$. Let $n$ be the least positive integer with
$nb=sa, \,1\leq n\leq l, \,1\leq s\leq k$ . Then $n\geq2$ and
$s\geq2$ by our assumptions. It is easy to see that

$$a, \,\ldots,  \,ka, \,\ldots, \,(k+[\frac{l}{n}]s)a,$$
$$b,\, a+b,\, \ldots, \, b+ka, \,\ldots, \,b+(k+[\frac{l-1}{n}]s)a,$$
$$\ldots \ldots$$
$$(n-1)b,\, \ldots, \,(n-1)b+ka, \,\ldots, \,(n-1)b+(k+[\frac{l-n+1}{n}]s)a$$
are all the distinct elements in $\sum(a^kb^l)$, and so
$$f(a^kb^l)=k+[\frac{l}{n}]s+1+k+[\frac{l-1}{n}]s+\cdots+1+k+[\frac{l-n+1}{n}]s$$
$$=n(k-s+1)+ls+s-1.$$
Since $n(k-s+1)+ls+s-1-2(k+l)=(n-2)(k-s)+(l-1)(s-2)+n-3$, we have
$f(a^kb^l)\geq 2(k+l)-1$ and the equality holds  if and only if
$n=s=2$, that is $2a=2b$. This completes the proof.

\end{proof}

{\bf Remark:} Note that if $a^kb^l\in \mathcal{A}^\star(G), \,k\geq
l\geq2$, then the conditions that $a^kb^l$ is smooth and $2a=2b$
cannot hold simultaneously. Otherwise, we may suppose that $2a=2b$
and $a^kb^l$ is $a$-smooth (the case that $a^kb^l$ is $b$-smooth is
similar), then $b=ta, \, 2\leq t\leq (k+1)$. It follows that
$b+(t-2)a=2(t-1)a=2b-2a=0, \,0<t-2\leq k-1$, which contradicts the
fact that $a^kb^l\in \mathcal{A}^\star(G)$.

\begin{lemma} {\rm [12, Lemma 2.9]}Let $S=a^kb^lg, \, k\geq l\geq1$ be a zero-sum free sequence over $G$ with
$b-a=g$ and  ${\rm ord}(g)=2$, then $f(S)=2(k+l)+1$.
\end{lemma}

\begin{lemma} Let $S_1\in \mathcal{F}(G)$ and $a, \,g\in G$ such
that $S=S_1a\in \mathcal{A}^\star(G)$, $S_1$ is $g$-smooth and $S$
is not $g$-smooth. Then $f(S)=2f(S_1)+1$.
\end{lemma}
\begin{proof} If $a\not\in<g>$, then
$\sum(S)=\sum(S_1)\cup\{a\}\cup(\sum(S_1)+a)$,  and so
$f(S)=2f(S_1)+1$.

If $a\in<g>$, we let $\sum(S_1)=\{g,\, \ldots, ng\}$, $a=tg,\,
t\in\mathbb{N}$, then $t\geq n+2$ by our assumptions. It follows
that $\sum(S)=\{g,\, \ldots, \, ng, \, tg, \, (t+1)g,\, \ldots, \,
(t+n)g\}$,  and so $f(S)=2f(S_1)+1$.\end{proof}

\begin{lemma} Let $k\geq 2$ be a positive integer and $a, \,b,\,c$
three distinct elements in $G$ such that
$a^kbc\in\mathcal{A}^\star(G)$ and $a^kbc$ is not $a$-smooth. Then
$f(a^kbc)\geq 2k+4$.
\end{lemma}
\begin{proof} Observe that $f(a^kbc)\geq
2k+4$ when $a^kbc$ is $b$ or $c$-smooth. We consider first the case
that $a^kb$ is $a$-smooth (the case that $a^kc$ is $a$-smooth is
similar). It is easy to see $f(a^kb)\geq k+2$, and so
$f(a^kbc)=2f(a^kb)+1\geq 2k+5$ by Lemma 2.7. Therefore we may assume
that both $a^kb$ and $a^kc$ are not $a$-smooth in the remaining
arguments. We divide the proof into three cases.

(i) If $a^k(b+c)$ is not $a$-smooth, then $a, \, \ldots, \,ka, \,
b,\, b+a,\, \ldots, \, b+c,\, b+c+a, \,\ldots, \,b+c+ka$ are
distinct elements in $\sum(a^kbc)$, and so
$$f(a^kbc)\geq k+k+1+k+1\geq 2k+4.$$

(ii) If neither $a^k(b-c)$ nor $a^k(c-b)$  is $a$-smooth, then $a,
\, \ldots, \,ka, \, b,\, b+a,\, \ldots, \, c,\, c+a, \,\ldots,
\,c+ka, \, b+c+ka$ are distinct elements in $\sum(a^kbc)$, and so
$$f(a^kbc)\geq k+k+1+k+1+1\geq 2k+5.$$

(iii) If $a^k(b+c)$ is  $a$-smooth and $a^k(b-c)$ (or $a^k(c-b)$  )
is $a$-smooth, then we have
$$ b+c=sa, \quad b-c=ta, \quad 1\leq s, \,t\leq k+1, \quad s\neq t.$$
It is easy to see that $a, \, \ldots, \,ka,\, (k+1)a, \,\ldots,
\,(k+s)a,$  \,$c, \,c+a, \, \ldots, \,c+(k+t)a$  are all distinct
elements in $\sum(a^kbc)$, and so
$$f(a^kbc)= k+s+k+t+1\geq 2k+4.$$
The second equality holds if and only if $(s, \,t)=(1,\,2)$ or $(2,
\,1)$. We are done.
\end{proof}
The following corollary follows immediately from Lemmas 2.1, and 2.7
and the proof of Lemma 2.9.
\begin{corollary} Let $k\geq 1$ be a positive integer and $a, \,b,\,c,\,d$
four distinct elements in $G$ such that
$a^kbcd\in\mathcal{A}^\star(G)$ and $a^kbcd$ is not $a$-smooth. Then
$f(a^kbcd)\geq 2k+6$.
\end{corollary}

\begin{lemma} Let $a, b, x$ be three distinct elements in $G$ such
that $a^kb^lx\in\mathcal{A}^\star(G), \, k\geq l\geq1$, $2a= 2b$,
and $x\neq a-b$, then $f(a^kb^lx)\geq 2(k+l+1)+1$.\end{lemma}
\begin{proof} If there are no distinct pairs $(m,\, n)\neq(0, \,0),
(m_1, \, n_1)\neq(0, \,0), \, 0\leq m, \, m_1\leq k, \,0\leq n, \,
n_1\leq l$ such that $ma+nb=m_1a+n_1b+x$, then
$\sum(a^kb^lx)=\sum(a^kb^l)\cup\{x\}\cup(\sum(a^kb^l)+x)$, and so
$f(a^kb^lx)=2f(a^kb^l)+1=4(k+l)-1\geq2(k+l+1)+1$.

If there are two distinct pairs $(m,\, n)\neq(0, \,0), (m_1, \,
n_1)\neq(0, \,0), \, 0\leq m, \, m_1\leq k, \,0\leq n, \, n_1\leq l$
such that $ma+nb=m_1a+n_1b+x$, then $x=a-b$ or $x=ua+b, \,1\leq
u\leq (k+l-1)$ or $x=vb, \, v\geq2$ or $x=ta, \,t\geq2$.

Let $x=ua+b, 1\leq u\leq (k+l-1)$, then $ a, \,\ldots, \,(k+l+u)a,
\, b, \, \cdots, \, b+(k+l+u)a$ are all distinct elements in
$\sum(a^kb^lx)$, and so $f(a^kb^lx)=2(k+l+u)+1\geq 2(k+l+1)+1$.

Let $x=vb, \, 2\leq v\leq (k+l)$ (the case that $x=ta, \,t\geq2$ is
similar). If $k$ is even, then $b, \, \ldots, \,(k+l+v)b, \, a, \,
a+b,\, \ldots, a+(k+l-2+v)b$ are all distinct elements in
$\sum(a^kb^lx)$, and so $f(a^kb^lx)=2(k+l+v-1)+1\geq 2(k+l+1)+1$. If
$k$ is odd, then $b, \, \ldots, \,(k+l+v-1)b, \, a, \, a+b,\,
\ldots, a+(k+l-1+v)b$ are all distinct elements in $\sum(a^kb^lx)$,
and so $f(a^kb^lx)=2(k+l+v-1)+1\geq 2(k+l+1)+1$. We are done.
\end{proof}

\begin{lemma} Let $a, b, x$ be three distinct elements in $G$ such
that $a^kb^2x\in\mathcal{A}^\star(G), \, k\geq 2$ and $a^kb^2x$ is
not $a$-smooth or $b$-smooth, then $f(a^kb^2x)=2k+5$ if and only if
$2a=2b$ and $x=b-a$.\end{lemma}
\begin{proof}We divide the proof into four cases.

{\bf Case 1} $a^kb^2$ is not smooth and $2b=sa, \,2\leq s\leq k$. If
$x=b-a$, then $ a, \, \ldots, \, (k+s)a, \, b-a, \, b, \, \ldots, \,
b+(k+s-1)a$ are all the distinct elements in $\sum(a^kb^2x)$, and so
$f(a^kb^2x)=2(k+s)+1$. If $x=ta, \,2\leq t\leq k$, then $ a, \,
\ldots, \, (k+s+t)a, \, b, \, \ldots, \, b+(k+t)a$ are all the
distinct elements in $\sum(a^kb^2x)$, and so
$f(a^kb^2x)=2(k+t)+s+1$.  If $x=ta+b, \,1\leq t\leq k$, then $ a, \,
\ldots, \, (k+s+t)a, \, b, \, \ldots, \, b+(k+t+s)a$ are all the
distinct elements in $\sum(a^kb^2x)$, and so
$f(a^kb^2x)=2(k+t+s)+1$. Therefore $f(a^kb^2x)=2k+5$ if and only if
$2a=2b$ and $x=b-a$ in this case.

{\bf Case 2} $a^kb^2$ is not smooth and $2b=sa, \, s> k$ or
$2b\not\in<a>$, then $f(a^kb^2)=3k+2$. If $k\geq3$, then
$f(a^kb^2x)\geq f(a^kb^2)+1=3k+3>2k+5$. If $k=2$ and $f(abx)=7$,
then $f(a^2b^2x)\geq f(abx)+f(ab)=7+3>2k+5$. If $k=2$ and $f(abx)=6$
(i.e., $x=a+b$ or $x=a-b$ or $x=b-a$), then it is easy to check that
$f(a^2b^2x)>2k+5$.

{\bf Case 3} $a^kb^2$ is smooth and $a^kb^2x$ is not smooth. If
$a^kb^2$ is  $a$-smooth, then
$f(a^kb^2x)=2f(a^kb^2)+1\geq2(k+2\times 2)+1> 2k+5$. If $a^kb^2$ is
$b$-smooth, then $f(a^kb^2x)=2f(a^kb^2)+1\geq2(2+2k)+1> 2k+5$.

{\bf Case 4} $a^kb^2x$ is $x$-smooth. We have
$f(a^kb^2x)\geq1+2k+2\times 3>2k+5$.

This completes the proof of the lemma.

\end{proof}

\section{Proofs of the Main Theorems}

To prove the main theorem of the present paper, we still need the
following two obviously facts on smooth sequences.

{\bf Fact 1} Let $r$ be a positive integer and $a\in G$. If $WT_i\in
\mathcal{A}^\star(G)$ is $a$-smooth for all $ i=1, \ldots, r$, then
$S=T_1\cdot \, \cdots \, \cdot T_rW$ is $a$-smooth.

{\bf Fact 2} Let $r, \, k, \, l$ be three positive integers and $a,
\, b$ two distinct elements in $G$. If $S\in \mathcal{A}^\star(G)$
is $a$-smooth and $a^kb^lT_i\in \mathcal{A}^\star(G)$ is $a$-smooth
or $b$-smooth for all $ i=1, \ldots, r$, then the sequence
$Sa^kb^lT_1\cdot \, \cdots \, \cdot T_r$ is $a$-smooth or
$b$-smooth.

{\bf Proof of Theorem 1.1:}

We start with the trivial case that $S = a^k$ with $k\in\mathbb{ N}$
and $a\in G$. Then $\sum(S) = \{a, \ldots, ka\}$, and since $S$ is
zero-sum free, it follows that $k < ord(a)$. Thus $S$ is $a$-smooth.

If $S=S_1g\in  \mathcal{A}^\star(G)$, where $g$ is an element of
order 2, then $f(S)\geq 2|S|-1$ by Lemma 2.4, and $f(S)\geq
f(S_1)+2$ since $\sum(S)\supseteq\sum(S_1)\cup\{g, \,
g+\sigma(S_1)\}$. If $S=S_1g_1g_2\in \mathcal{A}^\star(G)$, where
$g_1$ and $g_2$ are two elements of order 2, then  $f(S)\geq 2|S|$
since $\sum(S)\supseteq\sum(S_1g_1)\cup\{g_2, \, g_1+g_2, \,g_1+
g_2+\sigma(S_1)\}$. Therefore it suffices to determine all $S\in
\mathcal{A}^\star(G)$ such that $S$ does not contain any element of
order 2 and $f(S)\leq 2|S|-1$, and when $f(S)\leq 2|S|-1$, determine
all $Sg\in \mathcal{A}^\star(G)$ such that  $g$ is an element of
order 2 and $f(Sg)=2|S|+1$.

To begin with, we determine all $S\in \mathcal{A}^\star(G)$ such
that $S$ does not contain any element of order 2 and $f(S)\leq
2|S|-1$.  Let $S = a^xb^yc^zT$ with $x\ge y \ge z\ge \mathsf h(T)$
and  $a, \,b, c\not\in \supp(T)$. The case that $|\supp(S)|=2$
follows from Lemmas 2.3 and 2.5 and the remark after Lemma 2.5.
Therefore we may assume that $|\supp(S)|\geq3$ and $S$ does not
contain any element of order $2$ in the following arguments.

If $x=y=z$, then $S$ allows the product decomposition
$$S=S_1\cdot\,\cdots\,\cdot S_x,$$
where $S_i=abc\cdot\,\cdots, \,i=1, \,\ldots, \,x$ are squarefree of
length $|S_i|\geq3$. By Lemma 2.1, we obtain
$$f(S)\geq\sum_{i=1}^xf(S_i)\geq2\sum_{i=1}^x|S_i|=2|S|.$$

If $x\geq y>z\geq\mathsf h(T)$, or $x>y\geq z\geq\mathsf h(T)$, then
$S$ allows a product decomposition

\[
S = T_1 \cdot \,\cdots \,\cdot T_r W
\]
having the following properties:

\begin{itemize}

\item $r \geq 1$ and, for every $i \in [2, \,r]$, $S_i
\in \mathcal F (G)$
      is squarefree of length $|S_i| = 3$.

\item $W \in \mathcal F (G)$ has the form $W = a^k, \,k\geq 1$ or $W =
a^kb, \, k\geq1$  or $W=a^kb^l, \, k\geq l\geq 2$.

\end{itemize}
We choose  a product decomposition such that $k$ is the largest
integer in $W=a^k$ (or $a^kb$ or $a^kb^l, \, k\geq l\geq2$) among
all such product decompositions.  We divide the remaining proof into
three cases.

{\bf Case 1} $W=a^k, \,k\geq1$. If $T_i=xyz$ with $a\not\in\{x, \,
y, \,z\}$ for some $i, \,1\leq i\leq r$ such that $a^kxyz$ is not
$a$-smooth whenever $k>1$, then $S$ admits the product decomposition
$$
S = T_1 \cdot\, \cdots \,\cdot T_{i-1}T_i'T_{i+1} \cdot \,\cdots
\,\cdot T_r, $$ where $T_i, \,i=1, \,\ldots,  \,r$ have the
properties described above and $T_i'=a^{k}xyz$. By Lemma 2.1, and
Corollary 2.1, we get
$$f(S)\geq\sum_{j\neq i}^r f(T_j)+f(T_i')\geq\sum_{j\neq i}^r2|T_j|+2|T_i'|=2|S|.$$

If $T_i=axy$ for some $i, \,1\leq i\leq r$ such that $a^{k+1}xy$ is
not $a$-smooth, then $S$ admits the product decomposition
$$
S = T_1 \cdot\, \cdots \,\cdot T_{i-1}T_i'T_{i+1} \cdot \,\cdots
\,\cdot T_r, $$ where $T_i, \,i=1, \,\ldots,  \,r$ have the
properties described above and $T_i'=a^{k+1}xy$. By Lemmas 2.1 and
2.8, we get
$$f(S)\geq\sum_{j\neq i}^r f(T_j)+f(T_i')\geq\sum_{j\neq i}^r2|T_j|+2|T_i'|=2|S|.$$
Therefore we have proved that if $S$ is not $a$-smooth   and
$W=a^k$, then $f(S)\geq 2|S|$.

 {\bf Case 2} $W=a^kb, \,k\geq 1$.

  Let $T_i=xyz$ with $a\not\in\{x, \,
y, \,z\}$ for some $i, \,1\leq i\leq r$. If $k=1$, then
$T_iW=abxyz$. If $k=2$, then $T_iW=abx\cdot ayz$. If $k\geq3$ and
one sequence among three sequences  $a^{k-1}yz, \, a^{k-1}xz$, and
$a^{k-1}xy$, say, $a^{k-1}yz$ is not $a$-smooth, then $T_iW=abx\cdot
a^{k-1}yz$. It follows from Lemmas 2.1 and 2.8 that
$f(T_iW)\geq2|T_i|+2|W|$, and so $f(S)\geq2|S|$.

Let  $T_i=bxy$ for some $i, \,1\leq i\leq r$, then $k\geq2$. If
$k=2$, then $T_iW=abx\cdot aby$. If $k>2$ and $a^{k-1}by$ (or
$a^{k-1}bx$) is not $a$-smooth, then $T_iW=abx\cdot a^{k-1}by$ (or
$T_iW=aby\cdot a^{k-1}bx$). It follows from Lemmas 2.1 and 2.8 that
$f(T_iW)\geq2|T_i|+2|W|$, and so $f(S)\geq2|S|$.

Let  $T_i=abx$ for some $i, \,1\leq i\leq r$, then
$T_iW=a^{k+1}b^2x$. If $a^{k+1}b^2x$ is not $a$-smooth or
$b$-smooth, then by Lemma 2.10 we have $f(T_iW)\geq2|T_i|+2|W|$, and
so $f(S)\geq2|S|$.

Therefore we have proved that if $S$ is not $a$-smooth or
$b$-smooth, then $f(S)\geq 2|S|$ in this case.

{\bf Case 3} $W=a^kb^l, \,k\geq l\geq2$. If $2a\neq 2b$ and $a^kb^l$
is not smooth, then by Lemma 2.5 we have $f(W)\geq 2|W|$ and we are
done. Note that the conditions that  $2a=2b$ and $a^kb^l$ is smooth
cannot hold simultaneous. Here we omit the similar arguments as we
have done in Case 1.

{\bf Subcase 1}  $2a=2b$.

Let $T_i=xyz$ with $a\not\in\{x, \, y, \,z\}$ for some $i, \,1\leq
i\leq r$,  then $T_iW=abxy\cdot a^{k-1}b^{l-1}z$. It follows from
Lemmas 2.1 and 2.9 that $f(T_iW)\geq2|T_i|+2|W|$, and so
$f(S)\geq2|S|$.

Let  $T_i=byz$ for some $i, \,1\leq i\leq r$, then $k\geq l+1$,
$T_iW=aby\cdot a^{k-1}b^lz$.  It follows from Lemmas 2.1 and 2.9
that $f(T_iW)\geq2|T_i|+2|W|$, and so $f(S)\geq2|S|$.

Let  $T_i=abx$ for some $i, \,1\leq i\leq r$, then
$T_iW=a^{k+1}b^{l+1}x$. If $a^{k+1}b^2x$ is not $a$-smooth or
$b$-smooth, then by Lemma 2.10 we have $f(T_iW)\geq2|T_i|+2|W|$, and
so $f(S)\geq2|S|$.

{\bf Subcase 2} $a^kb^l$ is smooth, $a\neq 2b$, and $b\neq 2a$. Then
$W=(a^2b^2)^sW_1$, $W_1=a^{k_1}$ or $W_1=a^{k_1}b$. If
$S_1=SW^{-1}W_1$ is not $a$-smooth or $b$-smooth, then $f(S_1)\geq
2|S_1|$, and so by Lemmas 2.1 and 2.2 $f(S)\geq
sf(a^2b^2)+f(S_1)\geq8s+2|S_1|=2|S|$. If $S_1=SW^{-1}W_1$ is
$a$-smooth or $b$-smooth, then $S$ is $a$-smooth or $b$-smooth.

{\bf Subcase 3} $a=2b$.

Let $T_i=xyz$ with $a, \, b\not\in\{x, \, y, \,z\}$ for some $i,
\,1\leq i\leq r$,  then it is easy to see that
$f(T_iW)=f(a^kb^lxyz)=f(b^{2k+l}xyz)$. It follows from Corollary 2.1
that $b^{2k+l}xyz$ is $b$-smooth or $f(T_iW)\geq 2(|T_i|+|W|)$.

Let $T_i=bxy$ with $a, \,b\not\in\{x, \, y \}$ for some $i, \,1\leq
i\leq r$,  then  $f(T_iW)=f(a^kb^{l+1}xy)=f(b^{2k+l+1}xyz)$. It
follows from Lemma 2.8 that $b^{2k+l+1}xy$ is $b$-smooth or
$f(T_iW)\geq 2(|T_i|+|W|)$.

Let $T_i=abx$ with $a\neq x, \, b\neq x$ for some $i, \,1\leq i\leq
r$, then $f(T_iW)=f(a^{k+1}b^{l+1}x)=f(b^{2k+l+3}xyz)$. It follows
from Lemma 2.3 that $b^{2k+l+3}x$ is $b$-smooth or $f(T_iW)\geq
2(|T_i|+|W|)$.

{\bf Subcase 4} $b=2a$. Similar to Subcase 3.

Therefore we have proved that if $S$ is not $a$-smooth or
$b$-smooth, then $f(S)\geq 2|S|-1$ and $f(S)=2|S|-1$ if and only if
$S=a^kb$ or $S=a^kb^l, \, 2a=2b, \, k\geq l\geq2$.

Finally, when $f(S)\leq 2|S|-1$, we will  determine all $Sg\in
\mathcal{A}^\star(G)$ such that $g$ is an element of order 2 and
$f(Sg)=2|S|+1$.

(i) If $S$ is $a$-smooth (the case that $S$ is $b$-smooth is
similar), we set $\sum(S)=\{a, \, \ldots, \, na\}, \, n\leq 2|S|-1$,
then $g\not\in\sum(S)$ since $g$ is an element of order 2 and $Sg\in
\mathcal{A}^\star(G)$. It follows that
$\sum(Sg)=\sum(S)\cup\{g\}\cup\{g+\sum(S)\}$, and so $f(Sg)=2n+1$.
Therefore $f(Sg)\leq 2|S|+1$ if and only if $S=a^k$.

(ii) $S=a^kb$ is not smooth, by Lemma 2.8, $f(a^kbg)\leq 2k+1$ only
if $a^kbg$ is $a$-smooth, which is impossible since $g$ is an
element of order 2 and $a^kbg\in \mathcal{A}^\star(G)$.

(iii) $S=a^kb^l, \, 2a=2b, \, k\geq l\geq2$. The result follows from
Lemmas 2.5 and 2.9.

Therefore we have proved that if  $S=a^xb^y\cdot\,\cdots\in
\mathcal{A}^\star(G),  \,x\geq y\geq \cdots$, where $a, \,b, \,
\ldots$ are distinct elements of $G$ and $f(S)\leq 2|S|-1$, then $S$
is $a$-smooth or $b$-smooth or $S=a^kb, \, b\not\in\sum(a^k)$ or
$S=a^kb^l, \, k\geq l\geq2, 2a=2b$ or $S=a^kb^l, \, k\geq l\geq2,
2a=2b, \, g=a-b$. Theorem 1.1 is proved.

\hfill{$\Box$}

{\bf Acknowledgement:} The author wishes to thank Alfred Geroldinger
for sending the preprint \cite{al08} to him. He also thanks the
referee for his/her valuable suggestions.

$\; $
\noindent Pingzhi Yuan\\
School of Mathematics \\
South China Normal University\\
Guangdong, Guangzhou 510631\\
P.R.CHINA\\
e-mail:mcsypz@mail.sysu.edu.cn

\end{document}